\documentclass[psamsfonts,onesided,10pt]{amsart}

\usepackage{amssymb,amsmath,amsfonts}
\usepackage[all,arc,cmtip]{xy}
\usepackage{tikz}
\usepackage{enumerate}
\usepackage{mathrsfs}
\usepackage[margin=1in]{geometry}
\usepackage{algorithm}
\usepackage{graphicx}
\usepackage{algpseudocode}
\usepackage{booktabs}
\usepackage{longtable}
\usepackage{wrapfig}

\theoremstyle{definition}

\theoremstyle{remark}

\makeatletter
\let\c@equation\c@thm
\makeatother
\numberwithin{equation}{section}

\title{Algorithms for Solving Optimization Problems Arising from Deep Neural Net Models: Nonsmooth Problems}

\author{Vyacheslav Kungurtsev and Tomas Pevny}

\begin{document}

\begin{abstract}
Machine Learning models incorporating multiple layered learning networks have been seen to provide effective
models for various classification problems. The resulting optimization problem to solve for the optimal vector
minimizing the empirical risk is, however, highly nonconvex. This alone presents a challenge to application and development
of appropriate optimization algorithms for solving the problem. However, in addition, there are a number of interesting
problems for which the objective function is nonsmooth and nonseparable. In this paper, we summarize the primary challenges
involved, the state of the art, and present some numerical results on an interesting and representative class of problems.
\end{abstract}
\maketitle

\section{Introduction}
In recent years there has been a notable increase in the popularity of models in machine learning incorporating multiple 
layers of classifiers, referred to as Deep Neural Nets (DNNs). 
Appled to text classification, perception and identification, and a myriad of other settings, deep learning
has, after a prolonged slow start, showed impressive efficacy~\cite{lecun2015deep}. The models arising from 
deep learning, being the product of a layer of different models, result in a very nonconvex problem. In the age
of big data, the size of the problems that are desired to be solved is very large, necessitating the use
of stochastic or batch methods. 

Consider the general problem,
\[
\min_{x\in\mathbb{R}^n} F(x),
\]
where $F$ is generally nonconvex. In this paper in particular, we are interested in nonsmooth problems. Generally, the literature
on optimization algorithms for machine learning, even with deep learning in mind, has focused on smooth problems,
or, sometimes, smooth problems with an additional separable regularization term (usually the $l1$ norm)~\cite{bottou2016optimization}.
In the case of a simple regularizer, these are typically separable, that is, decomposable as a sum of functions
of individual variables or blocks (e.g., $\|x\|_1 = \sum_{i=1}^N |x_1|$), which permits a closed form 
solution of a \textbf{prox} operator, as in the iterative shrinkage algorithm. Coordinate descent type algorithms
have proven successful for problems with nonconvex loss functions and separable 
regularizers~\cite{richtarik2016parallel, daneshmand2015hybrid}. Although extensions to include stochastic or batch
subproblem/step formulation are necessary to truly apply them to huge scale deep learning problems, beyond that
they can only converge to a coordinate-wise stationary point, i.e., a point $x^*$ such that there exists $\xi\in\partial_i F(x^*)$
with $(x_i-x_i^*)^T \xi\ge 0$ for all $x$, $i$, but this does not necessarily mean stationarity, wherein there exists $\xi\in\partial F(x^*)$
with $(x-x^*)^T\xi\ge 0$ for all $x$.

However, in many applications, the actual loss function itself can be nonsmooth. In particular, some upper layers
can involve an activation function, defined as involving terms of the form $\max(0,a^T x+b)$. This presents
a large scale optimization problem that is nonconvex, nonsmooth, and nonseparable, a particularly challenging
class of problems with no algorithm that stands out as particularly suitable for reliably and efficiently 
solving them. 

\section{Nonsmooth, nonconvex problems, state of the art and prognosis}
In general, solving nonconvex nonsmooth problems is a formidable task. Provably convergent algorithms for solving such
problems are uncommon. One of the most robust such algorithms, gradient sampling~\cite{burke2005robust} requires $2n$ samples of 
the gradient to be taken at each iteration, which implies that its performance scales very poorly with the size of the problem. 
Yet this algorithm is regularly used in areas of robust control, for instance. Aside from this, the provably convergent
algorithms include classes of bundle methods and some other subgradient-approximation type algorithms, namely discrete gradient
and pseudosecant methods. Thorough numerical experiments given in~\cite{bagirov2014introduction} suggest that for large
scale problems, the limited memory bundle method~\cite{haarala2004new} and the quasisecant method~\cite{bagirov2013subgradient}
were the most robust and efficient. 

However, the quasisecant method requires the computation of quasisecants at each iteration, which is a very problem
dependent and potentially complicated procedure~\cite{bagirov2010quasisecant}. Given the convoluted structure
of deep neural net models, the computation of quasisecants at each iteration of an optimization for solving
such a problem would be intractable. With respect to bundle methods, each iteration involves solving a quadratic program (QP)
that summarizes the function locally from the historical information of (sub)gradients. Often in applications
involving big data, computation of an entire gradient vector is impossible, given the storage difficulties associated
with all of the data batch functions. Thus either stochastic, or batch averages, or gradient approximations must be
used. The robustness of using such approximations in a bundle method can be questionable. Although there are 
bundle methods incorporating inexact information~\cite{de2014convex}, experience with applying them to this class
of problems is limited, at best. We shall see in the next Section~\ref{s:num} that it is, at best, not straightforward.

\section{Numerical Results}\label{s:num}
Recall that by Rademacher's Theorem, it holds that for every point outside a set of measure zero, a unique
gradient vector can be computed~\cite{bessis1999partial}. This implies that gradient methods, or more specifically,
stochastic gradient and batch gradient methods are, in practice, well-defined for nonsmooth problems. That said, there
exists no convergence theory for the basic variation of the (sub)gradient method, even with a diminishing step-size,
for nonsmooth nonconvex problems. In terms of practical implementation, however, it is practical to suggest them
as a method to try for this class of problems, in their generic form, as they are the standard bearer for large
scale optimization for problems arising in machine learning in general. 

Just as basic SGD methods have been seen to converge slowly for highly nonconvex problems with a lot
of saddle points, however, we expect that their practical convergence speed for highly nonconvex, and thus
ill-conditioned problems, to be fairly slow. Thus, just as in the smooth case, we expect that there is a case
for higher-order algorithms, but with a number of difficulties in reliable and efficient implementation. 
Note that, we at least at least one application area of strong interest, with the activation function mentioned above,
to be piecewise linear. This implies that the second derivative is zero almost everywhere (and then undefined on
a set of measure zero), and thus Newton methods (or stochastic variations thereof) are entirely useless. 
However, despite a lack of convergence theory, quasi-Newton, specifically BFGS methods have numerically demonstrated
strong efficacy for solving general nonsmooth, nonconvex problems~\cite{lewis2013nonsmooth}. Given the scale
of the problems of interest, storing dense Hessians is impractical, and thus the limited memory variety of BFGS~\cite{liu1989limited}
has been tested on optimizations problems arising from DNN models to good success~\cite{dean2012large, ngiam2011optimization}.

We compared the performance of ADAM~\cite{kingma2014adam}, a contemporary stochastic gradient algorithm that scales
the direction based on accumulated history of apparent conditioning, sum of functions optimizer (SFO)~\cite{sohl2014fast}, a second order
BFGS-type algorithm using batches, and the limited memory bundle method (LMBM)~\cite{haarala2004new} on a model of 
deep learning of security data from Cisco in Prague, in the Czech Republic. We studied their performance, 
indicating the scalability, on a variety of sizes of batches, each composed of several hundred samples. 
The LMBM was run in an entirely sequential and undistributed manner, i.e., treating the problem in its
entirety, while the other two methods are stochastic or single-batch based. 

\begin{figure}
\caption{\label{fig-valplot} Plot comparing the final value after 1000 iterations of ADAM, SFO, and LMBM.}
\begin{center}
\includegraphics[scale=0.52]{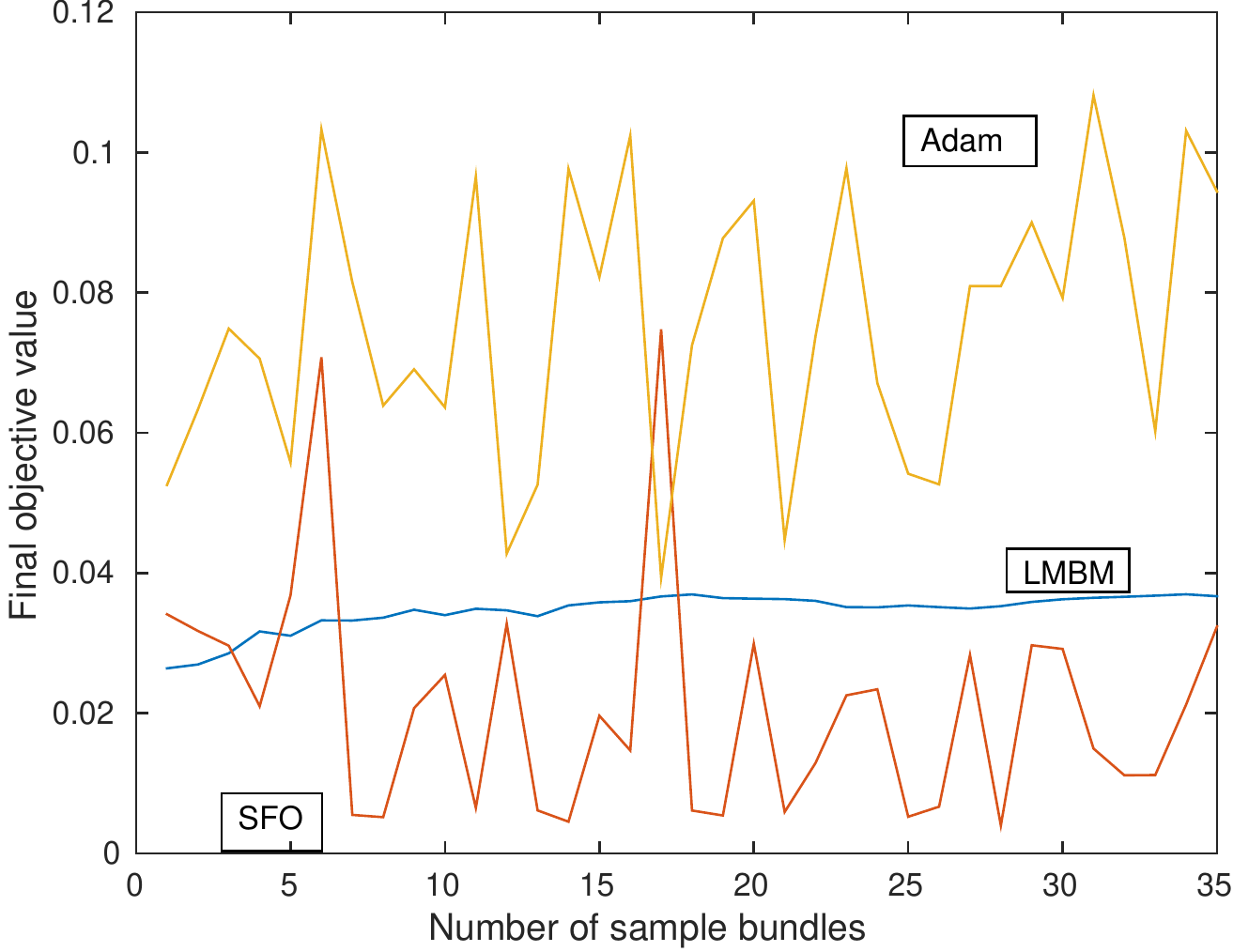}
\end{center}
\end{figure}

We see in Figure~\ref{fig-valplot} the final value of each algorithm after 1000 iterations. Overall it is clear
that the SFO algorithm results in the best objective value, LMBM is second overall, but more consistent/robust, 
and finally ADAM is last.

\begin{figure}
\caption{\label{fig-timeplot} Plot comparing the time to convergence of ADAM, SFO, and LMBM.}
\begin{center}
\begin{tabular}{c c} 
\includegraphics[scale=0.35]{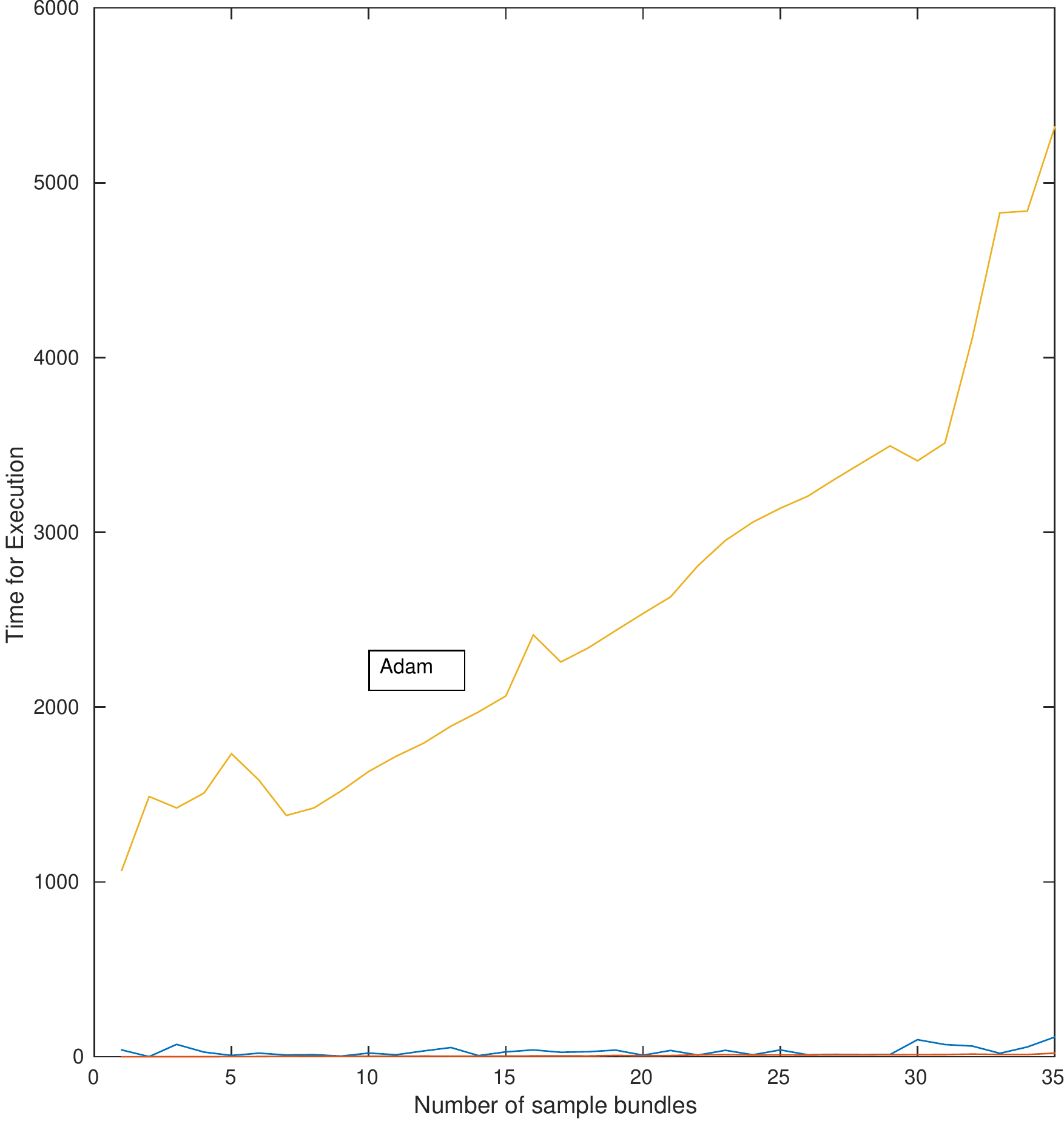}
&
\includegraphics[scale=0.51]{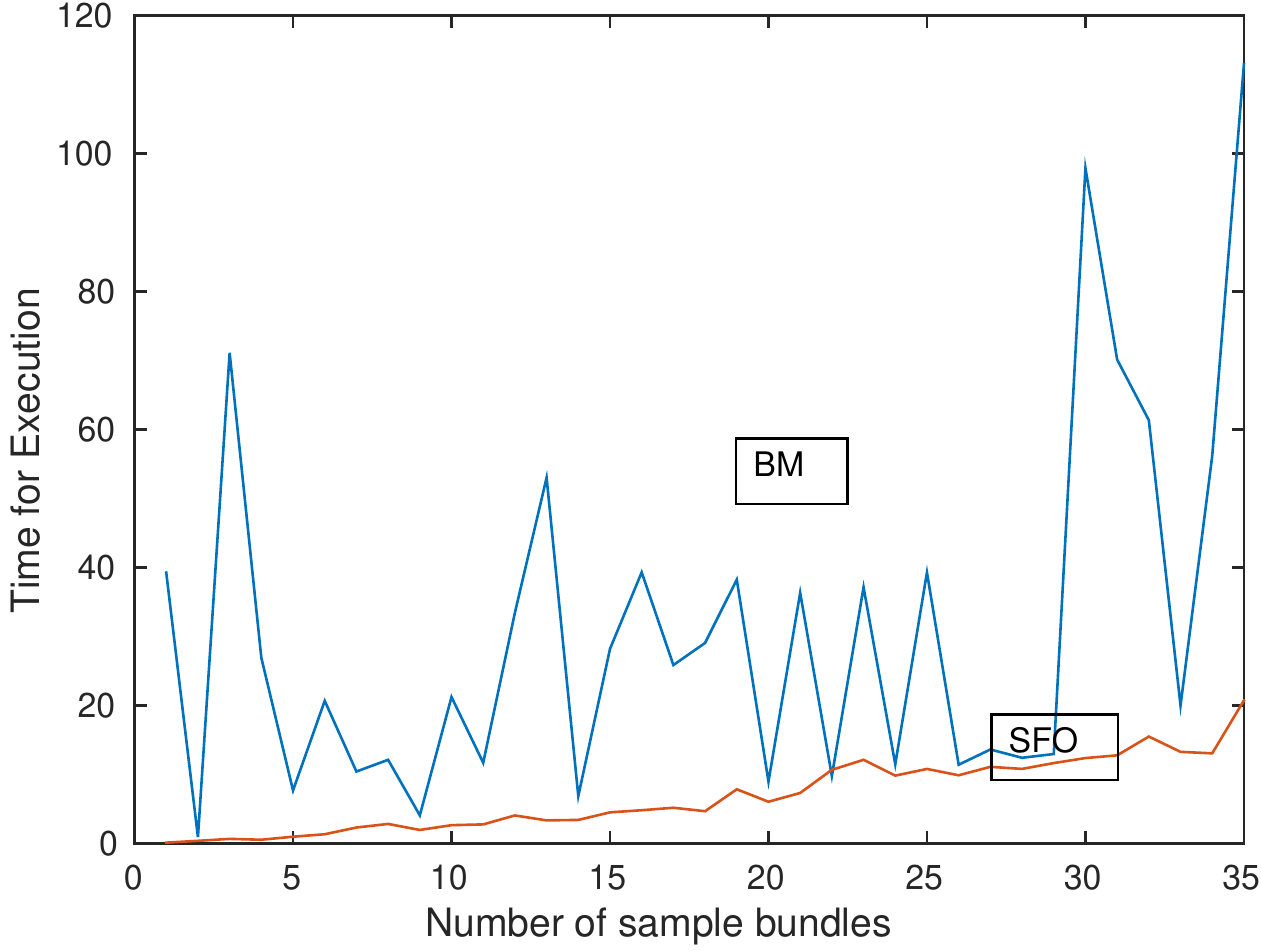}
\end{tabular}
\end{center}
\end{figure}

Taking the final value that ADAM generates after 1000 iterations, we now plot the time until each algorithm
reaches this value in Figure~\ref{fig-timeplot}. This also suggests that SFO is the most efficient algorithm,
followed by LMBM and ADAM as a distant second.

After 35 batches, the required computation and storage that was needed for the bundle method was too much.
We attempted to modify the LMBM to handle larger problem sizes. Thus far, none of the attempts we implemented, including,
\begin{enumerate}
\item computing the entire gradient by computing each batch gradient and averaging,
\item choosing a gradient from a random batch each iteration,
\item running LMBM to minimize the problem for each batch sequentially,
\end{enumerate}
resulted in an efficiently and quickly converging algorithm. This does not necessarily suggest that
huge-scale bundle methods are impossible, but their implementation is certainly not straightforward. 

\section{Discussion}
A provably convergent algorithm for huge scale nonconvex, nonsmooth, and nonseparable optimization problems
is, for the most part, a open problem. The demands from applications in machine learning, with the arrival
of the widespread use of deep neural nets, has shifted the necessity of solving large convex problems
to nonconvex ones. The optimization literature has matched this demand to some degree, but the added
challenge of solving nonsmooth problems, aside from simple separable regularizers, appears to be far more
difficult. Numerical experience suggests the efficacy of carefully designed quasi-Newton methods,
but convergence is not theoretically guaranteed. The top competing algorithms for nonconvex, nonsmooth
optimizations on a (smaller, but still) large scale, quasi-secant and bundle methods, are not straightforward
to extend to a larger scope of problem size, in the first case due to the challenge of computing the 
step, and in the second due to memory requirements. We foresee that a proximal type bundle method with 
stochastic information can be developed in the future as a provably convergent algorithm for these
cases of problems, and we will focus on this as future work, even if thus far such a procedure 
remains elusive.

\bibliographystyle{plain}
\bibliography{refs}
\end{document}